\documentclass[12pt]{article} 

\usepackage{amsmath}
\usepackage{latexsym}
\usepackage{url}
\usepackage{color}
\numberwithin{equation}{section}

\def\R{{\bf R}}

\def\N{{\bf N}}

\def\d{\displaystyle}
\def\e{{\varepsilon}}

\def\wt{\widetilde}

\def\p{\partial}

\newtheorem{thm}{Theorem}[section]

\newtheorem{rem}{Remark}[section]

\title{Blow-up of classical solutions of quasilinear wave equations in one space dimension}
\author{
Yuki Haruyama
\footnote{
Master course, Mathematical Institute,
Tohoku University,
Aoba, Sendai 980-8578, Japan.
email: yuki.haruyama.t6@dc.tohoku.ac.jp}
\ and
Hiroyuki Takamura
\footnote{Mathematical Institute,
Tohoku University,
Aoba, Sendai 980-8578, Japan.
e-mail: hiroyuki.takamura.a1@tohoku.ac.jp.}
}
\date{
\[
\begin{array}{ll}
\mbox{\footnotesize{\bf Keywords:}}
& \mbox{\footnotesize quasilinear wave equation, one dimension, blow-up, classical solution, lifespan}\\
\mbox{\footnotesize{\bf MSC2020:}}
& \mbox{\footnotesize primary 35L72, secondary 35B44}\\
\end{array}
\]
}

\begin{document}
\maketitle

%%%%%%%%%%%%%%%%%%%%%%%%%%%%%%%%%%%%%%%%%%%%%%%
%%%%%%%%%%%%%%%%%%% ABSTRACT %%%%%%%%%%%%%%%%%%%%%
%%%%%%%%%%%%%%%%%%%%%%%%%%%%%%%%%%%%%%%%%%%%%%% 

\begin{abstract}
This paper studies the upper bound of the lifespan of classical solutions
of the initial value problems for one dimensional wave equations
with quasilinear terms of space-, or time-derivatives of the unknown function.
The result for the space-derivative case guarantees
the optimality of the general theory for nonlinear wave equations,
and its proof is carried out by combination of ordinary differential inequality
and iteration method on the lower bound of the weighted functional of the solution.
\end{abstract}

%%%%%%%%%%%%%%%%%%%%%%%%%%%%%%%%%%%%%%%%%%%%%%%%%%
%%%%%%%%%%%%%%%%%%%%% SECTION1 %%%%%%%%%%%%%%%%%%%%%%%
%%%%%%%%%%%%%%%%%%%%%%%%%%%%%%%%%%%%%%%%%%%%%%%%%%

\section{Introduction}
In this paper, we first consider the following initial value problems of a scaler unknown function
$v=v(x,t)$ of spatial variables;
\begin{equation}
\label{IVP}
\left\{
\begin{array}{ll}
v_{tt} - v_{xx} = A|v_x|^{p-2}v_xv_{xx}
&\mbox{in}\ \in\R\times(0,T),\\
v(x,0)= \e f(x),\ v_t(x,0)= \e g(x),
& x\in\R,
\end{array}
\right.
\end{equation}
where $p>1, A,T>0$.
We assume that $f$ and $g$ are given smooth functions of compact support
and a parameter $\e>0$ is \lq\lq small enough".
We are interested in the lifespan $T(\e)$, the maximal existence time,
of classical solutions of (\ref{IVP}).
The significant meaning to consider this problem is to show the optimality
of the general theory for nonlinear wave equation on the lifespan estimates
especially for quasilinear case with any power nonlinearity.
See Takamura \cite{Takamura} for the most recent results of the theory.
\par
More precisely, the trigger of our motivation to investigate (\ref{IVP}) comes from
Sasaki, Takamatsu and Takamura \cite{STT23} in which the equation
\begin{equation}
\label{eq:STT}
v_{tt}-v_{xx}=|v_x|^p
\end{equation}
is considered and its result is that the lifespan $T(\e)$ has to satisfy
\begin{equation}
\label{est:STT}
C_1\e^{-(p-1)}\le T(\e)\le C_2\e^{-(p-1)},
\end{equation}
where $C_1$ and $C_2$ are positive constants independent of $\e$.
We remark  that the proof of the lower bound in (\ref{est:STT}) has a trivial error.
One may read it correctly by replacing $w_t$ by $w_x$ in the function space and
$\overline{L'}$ by $L'$ in the convergence of the sequence $\{w_j\}$.

\par
We know that (\ref{est:STT}) is also available for the equation
\begin{equation}
\label{eq:Zhou}
v_{tt}-v_{xx}=|v_t|^p.
\end{equation}
Such a result  is due to Zhou \cite{Zhou01} for the upper bound
and Kitamura, Morisawa and Takamura \cite{KMT23} for the lower bound.
We note that (\ref{eq:Zhou}) is important to show the optimality
of the general theory for nonlinear wave equations by
Li, Yu and Zhou \cite{LYZ91, LYZ92}.
Recently, the general theory is improved by \lq\lq combined effect".
See Takamatsu \cite{Takamatsu} for this direction
as well as Morisawa, Sasaki and Takamura \cite{MST23, MST23e},
Kido, Takamatsu, Sasaki and Takamura \cite{KSTT} for its optimality.
Or, see Takamura \cite{Takamura} as a review on this research field.
The advantage to investigate (\ref{eq:STT}) beyond (\ref{eq:Zhou}) is not only to cover
the optimality of the theory by model equations as many as possible,
but also to extend the analysis on the \lq\lq blow-up boundary".
See Sasaki \cite{Sasaki18} for details.
The key fact is that a point-wise positiveness of the solution is available for (\ref{eq:Zhou}),
but (\ref{eq:STT}) is out of the case.
In this way, the semilinear model equations are well-studied.

\par
On contrary, as far as the authors know,
the quasilinear model in one space dimension has been studied only for the equation
\begin{equation}
\label{eq:John}
v_{tt}-c^2(v_x)v_{xx}=0
\end{equation}
by John \cite{John81}, or John \cite{John_book} on page 7,
where a given smooth function $c$ satisfies
\begin{equation}
\label{c_John}
c(v_x)>0\quad\mbox{and}\quad c'(v_x)\neq0.
\end{equation}
In this situation, the method of the characteristic curves gives us the blow-up of the solution
and the lifespan estimate of
\begin{equation}
\label{est:John}
T(\e)\le C\e^{-1},
\end{equation}
where $C$ is a positive constant independent of $\e$.
We note that the setting $c(v_x)=\sqrt{1+A|v_x|^{p-2}v_x}$
in our problem (\ref{IVP}) cannot be adopted to the assumption (\ref{c_John}) except for $p=2$.
Related to this exceptional case,
the equation of a different type,
\[
v_{tt}-\Delta v=2v_tv_{tt},
\]
in three space dimension is studied by John \cite{John_book} on page 27,
and it is easy to reduce its result and proof to one dimensional case.
We will see it at the last section of this paper by studying generalized equations like the one in (\ref{IVP}),
\begin{equation}
\label{eq:John_t}
v_{tt}-v_{xx}=A|v_t|^{p-2}v_tv_{tt}.
\end{equation}
We finally remark that another type of the quasilinear model,
\[
u_{tt}=c(u)^2u_{xx}+\lambda c(u)c'(u)u_x^2
\]
with a parameter $\lambda\in(0,1]$, is studied by Sugiyama \cite{Sugiyama2022}.
Under the uniform positive assumption on $c(0)$ and $c'(0)$,
he obtains that a blow-up time of the solution is bounded and independent of $\e$
due to speciality of the data with a concentration at the origin,
which is different from our setting on the initial data.
See the introduction of \cite{Sugiyama2022} for its details and references therein.

\par
Our main purpose in this paper is to extend John's result on (\ref{eq:John}) to (\ref{IVP}),
and to obtain the same upper bound in (\ref{est:STT}).
The optimality of such a result is already guaranteed by the general theory
when $p$ is any even integer.
In fact, the theory provides us the fact that
the equation, $v_{tt}-v_{xx}=Av_x^{p-1}v_{xx}$ with $p=2,4,6,\ldots$,
has the same lower bound of the lifespan as in (\ref{est:STT}).
See \cite{LYZ91, LYZ92} or \cite{Takamura} for this fact.
But, for other values of $p$ including fractional values,
the optimality of our result is still open
because the general theory is available only for smooth nonlinear terms.
We note that our method is different from \cite{John81, John_book}
after a reduction to a wave equation with positive terms
as we employ the argument in Sasaki, Takamatsu and Takamura
\cite{STT23} which is based on the weighted functional method
by Rammaha \cite{Rammaha95, Rammaha97}.
We also note that the analysis including the comparison argument
on the reduced ordinary differential inequality in \cite{Rammaha95, Rammaha97}
is replaced with a simple iteration argument.
Such a hybrid method was first introduced by Lai and Takamura \cite{LT18}
to avoid any restriction on the lower bound of the exponent in power-type nonlinear terms
in showing a blow-up result for semiliner damped wave equations.

\par
This paper is organized as follows.
The statements of our results are described in Section 2.
Section 3 is devoted to the proof of the upper bound of the lifespan of the classical solution to (\ref{IVP})
while Section 4 is devoted to the one in the case where the equation is (\ref{eq:John_t}).
Finally, the concluding remarks and future problems are given in Section 5.

%%%%%%%%%%%%%%%%%%%%%%%%%%%%%%%%%%%%%%%%%%%%%%%%%%
%%%%%%%%%%%%%%%%%%%%% SECTION2 %%%%%%%%%%%%%%%%%%%%%%%
%%%%%%%%%%%%%%%%%%%%%%%%%%%%%%%%%%%%%%%%%%%%%%%%%%

\section{Theorems}
All the results expected to be shown are listed in this section.
In any case, our assumption on the initial data is $(f,g)\in C_0^2(\R)\times C_0^1(\R)$ satisfying
\begin{equation}
\label{supp_data}
\mbox{supp}\ (f,g)\subset\{x\in\R:|x|\le \sigma\}\ (\sigma\ge1).
\end{equation}

\par
Our main result is the following which will be proved in the next section.

\begin{thm}
\label{thm:main}
Assume (\ref{supp_data}) and that $f$ and $g$ satisfy
\begin{equation}
\label{positive1}
f(x), g(x)\geq0,\ \mbox{and}\ f(x)\not\equiv0,
\end{equation}
or
\begin{equation}
\label{positive2}
f(x), g(x)\geq0,\ \mbox{and}\ g(x)\not\equiv0.
\end{equation}
Then, there exists a positive constant $\e_0=\e_0(f,g,p,A,\sigma)>0$ such that
any classical solution of (\ref{IVP}) in the time interval $[0,T]$ cannot exist
as far as $T$ satisfies
\begin{equation}
\label{upper-bound}
T>C\e^{-(p-1)},
\end{equation}
where $0<\e\leq\e_0$, and $C$ is a positive constant independent of $\e$.
\end{thm}

\par
As stated in Introduction, we have a similar result on the equation of time-derivatives
whose proof will be given at the last section.

\begin{thm}
\label{thm:t-derivative}
Even if the equation in (\ref{IVP}) is replaced with (\ref{eq:John_t}),
the conclusion of Theorem \ref{thm:main} is still valid.
\end{thm}

\begin{rem}
The estimate (\ref{upper-bound}) implies that
\[
T (\e)\le C\e^{-(p-1)}\quad \mbox{for}\ 0 < \e\le\e_0.
\]
\end{rem}

%%%%%%%%%%%%%%%%%%%%%%%%%%%%%%%%%%%%%%%%%%%%%%%%%%
%%%%%%%%%%%%%%%%%%%%% SECTION3 %%%%%%%%%%%%%%%%%%%%%%%
%%%%%%%%%%%%%%%%%%%%%%%%%%%%%%%%%%%%%%%%%%%%%%%%%%

\section{Proof of Theorem \ref{thm:main}}
First of all, we note that we may assume $A=p$ without loss of generality
by scaling of the unknown function.
We shall employ the contradiction argument.
Let $v$ be a classical solution of (\ref{IVP}) in the time interval $[0,T)$.
It is well-known that
\begin{equation}
\label{supp_v}
 \mbox{supp}\ v(x,t)\subset\{(x,t)\in\R\times[0,T):|x|\le t+\sigma\}.
\end{equation}
See Appendix in John \cite{John_book} for example.
Hence
\begin{equation}
\label{u1}
u(x,t):=\int_{-\infty}^xv(y,t)dy
\end{equation}
is well-defined in $\R\times[0,T)$ and
\begin{equation}
\label{supp_u1}
 \mbox{supp}\ u(x,t)\subset\{(x,t)\in\R\times[0,T):x\ge-t-\sigma\}.
\end{equation}
Then, it follows from the equation in (\ref{IVP}) that
\[
\frac{\partial}{\partial x} (u_{tt}-u_{xx}-|u_{xx}|^p)=0
\quad\mbox{in}\ \R\times(0,T).
\]
Integrating this identity in $x$ and making use of (\ref{supp_u1}), we obtain that 
\[
\left\{
\begin{array}{ll}
u_{tt}-u_{xx}=|u_{xx}|^p  &\mbox{in}\ \R\times(0,T),\\
\d u(x,0)=\e\int_{-\sigma}^{x}f(y)dy,\  u_t(x,0)=\e\int_{-\sigma}^{x}g(y)dy, & x \in \R.
\end{array}
\right.
\]
The function $u$ is a classical solution of this problem, so that $u$ has to satisfy an integral equation,
\begin{equation}
\label{eq:integral}
u=u_0+L(|u_{xx}|^p)\quad\mbox{in}\ \R\times[0,T),
\end{equation}
where we set
\begin{equation}
\label{u_0L}
\begin{array}{rl}
u_0(x,t)
&\d:= \frac{1}{2} \left\{u(x+t,0)+u(x-t,0)+\int_{x-t}^{x+t}u_t(y,0)dy\right\},\\
L(w)(x,t)
&\d:=\frac{1}{2}\int_{0}^{t} ds \int_{x-t+s}^{x+t-s}w(y,s)dy.
\end{array}
\end{equation}

\par
Now, following Rammaha \cite{Rammaha95, Rammaha97},
we shall introduce a new function $H=H(t)$ defined by
\begin{equation}
\label{H}
H(t):= \int_{0}^{t} (t-s)ds\int_{-s-\sigma}^{-s-\sigma_0}u(x,s) dx,
\end{equation}
where $\sigma_0$ is a constant satisfying $0<\sigma_0<\sigma$ which will be fixed appropriately later.
Then, it follows that
\[
H'(t)= \int_{0}^{t} ds \int_{-s-\sigma}^{-s-\sigma_0}u(x,s)dx,
\quad
 H''(t)= \int_{-t-\sigma}^{-t-\sigma_0} u(x,s) dx.
 \]
 Hence (\ref{eq:integral}) yields that
\begin{equation}
\label{H''}
\begin{array}{ll}
H''(t)=
&\d\frac{1}{2}\int_{-t-\sigma}^{-t-\sigma_0}\left\{u(x+t,0)+u(x-t,0)+\int_{x-t}^{x+t}u_t(y,0)dy\right\}dx\\
&\d+\frac{1}{2}F(t),
\end{array}
\end{equation}
where $F=F(t)$ is defined by
\begin{equation}
\label{F}
F(t):=\int_{-t-\sigma}^{-t-\sigma_0} dx\int_{0}^{t}ds\int_{x-t+s}^{x+t-s}|u_{xx}(y,s)|^pdy.
\end{equation}

\par
First, we shall prove the theorem under the assumption (\ref{positive1}).
Then, we have that
\[
H''(t)\ge\frac{1}{2} \int_{-t-\sigma}^{-t-\sigma_0}u(x+t,0)dx
= \frac{1}{2}\int_{-\sigma}^{-\sigma_0}u(x,0)dx
= \frac{\e}{2} \int_{-\sigma}^{-\sigma_0}dx\int_{-\sigma}^{x}f(y)dy
\]
which implies that
\[
H''(t)\ge\frac{\e}{2}\int_{-(\sigma+\sigma_0)/2}^{-\sigma_0}dx
\int_{-\sigma}^{-(\sigma+\sigma_0)/2}f(y)dy=C_f\e\quad\mbox{for}\ t\in[0,T),
\]
where $\sigma_0$ is fixed as
\begin{equation}
\label{C_f}
C_f:=\frac{\sigma-\sigma_0}{4}\int_{-\sigma}^{-(\sigma+\sigma_0)/2} f(y) dy>0.
\end{equation}
We note that we may assume $f\not\equiv0$ in $(-\sigma,-(\sigma+\sigma_0)/2)$ without loss of generality
by parallel translation in $x$.
Therefore it follows from $H(0)=H'(0)=0$ that
\begin{equation}
\label{est:H}
H(t)\ge Bt^2\quad\mbox{for}\ t\in[0,T),
\end{equation}
where $B:=C_f\e/2$.
This will be the first step of the iteration procedure.

\par
Next, we shall derive an iteration frame from (\ref{H''}).
Rewrite $F$ in (\ref{F}) as
\[
F(t)=\int_{0}^{t}ds\int_{-t-\sigma}^{-t-\sigma_0}dx\int_{x-t+s}^{x+t-s}|u_{xx}(y,s)|^pdy
\]
and assume that
\begin{equation}
\label{sigma_1}
t\ge\sigma_1:= \frac{\sigma-\sigma_0}{2}.
\end{equation}
Then, switching the order of $x$-$y$ integral, we have that
\[
\begin{array}{l}
\d\int_{-t-\sigma}^{-t-\sigma_0} dx\int_{x-t+s}^{x+t-s} |u_{xx}(y,s)|^p dy\\
\d=\left(\int_{-2t+s-\sigma}^{-2t+s-\sigma_0}dy\int_{-t-\sigma}^{y+t-s}dx
+\int_{-2t+s-\sigma_0}^{-s-\sigma}dy\int_{-t-\sigma}^{-t-\sigma_0}dx\right)|u_{xx}(y,s)|^p\\
\d\quad+\int_{-s-\sigma}^{-s-\sigma_0}dy\int_{y-t+s}^{-t-\sigma_0}|u_{xx}(y,s)|^pdx\\
\d\ge\int_{-s-\sigma}^{-s-\sigma_0}(-s-\sigma_0-y)|u_{xx}(y,s)|^pdy
\end{array}
\]
for $0\leq s \leq t-\sigma_1 $.
Similarly, we also have that
\[
\begin{array}{l}
\d\int_{-t-\sigma}^{-t-\sigma_0}dx\int_{x-t+s}^{x+t-s} |u_{xx}(y,s)|^p dy\\
\d=\int_{-2t+s-\sigma}^{-s-\sigma}dy\int_{-t-\sigma}^{y+t-s}|u_{xx}(y,s)|^pdx\\
\d\quad+\left(\int_{-s-\sigma}^{-2t+s-\sigma_0}dy\int_{y-t+s}^{y+t-s}dx
+\int_{-2t+s-\sigma_0}^{-s-\sigma_0}dy\int_{y-t+s}^{-t-\sigma_0}dx\right)|u_{xx}(y,s)|^p\\
\d\ge\int_{-s-\sigma}^{-2t+s-\sigma_0}2(t-s)|u_{xx}(y,s)|^pdy
+\int_{-2t+s-\sigma_0}^{-s-\sigma_0}(-s-\sigma_0-y)|u_{xx}(y,s)|^pdy            
\end{array}
\]
for  $t-\sigma_1 \leq s \leq t$.
Summing up the estimates above and making use of
\[
1= \frac{-s-\sigma_0-y}{-s-\sigma_0-y}
\ge\frac{-s-\sigma_0-y}{\sigma-\sigma_0}
\ge\frac{-s-\sigma_0-y}{2t}
\quad\mbox{for}\ y\ge-s-\sigma
\]
in supp $u(y,s)$ by (\ref{supp_u1}), we obtain that
\[
\begin{array}{ll}
F(t)\ge&
\d\int_{0}^{t-\sigma_1}ds\int_{-s-\sigma}^{-s-\sigma_0}(-s-\sigma_0-y)|u_{xx}(y,s)|^pdy\\
&\d+ \int_{t-\sigma_1}^{t}\frac{t-s}{t}ds \int_{-s-\sigma}^{-2t+s-\sigma_0}(-s-\sigma_0-y)|u_{xx}(y,s)|^pdy\\
&\d+ \int_{t-\sigma_1}^{t}ds\int_{-2t+s-\sigma_0}^{-s-\sigma_0}(-s-\sigma_0-y)|u_{xx}(y,s)|^pdy
\end{array}
\]
which yields that
\[
F(t)\ge\frac{1}{t}\int_{0}^t(t-s)ds\int_{-s-\sigma}^{-s-\sigma_0}(-s-\sigma_0-y)|u_{xx}(y,s)|^pdy
\quad\mbox{for}\ t\ge\sigma_1. 
\]
Therefore, neglecting the first integral in (\ref{H''}), we obtain that
\begin{equation}
\label{H''bound}
H''(t)\ge\frac{1}{2t} \int_{0}^{t} (t-s)ds \int_{-s-\sigma}^{-s-\sigma_0} (-s-\sigma_0-y)|u_{xx}(y,s)|^pdy
\quad \mbox{for}\ t\ge\sigma_1.
\end{equation}

\par
Now we shall estimate the right-hand side of (\ref{H''bound}) by $H$.
It follows from the integration by parts and (\ref{supp_u1}) that
\[
\begin{array}{ll}
H(t)
&\d=-\int_{0}^{t} (t-s)ds\int_{-s-\sigma}^{-s-\sigma_0}\p_y(-s-\sigma_0-y) u(y,s)dy\\
&\d=\int_{0}^{t} (t-s)ds\int_{-s-\sigma}^{-s-\sigma_0} (-s-\sigma_0-y)u_x(y,s)dy\\
&\d=-\int_{0}^{t}(t-s)ds \int_{-s-\sigma}^{-s-\sigma_0} \p_y\left\{\frac{1}{2}(-s-\sigma_0-y)^2\right\}u_x(y,s)dy\\
&\d=\int_{0}^{t} (t-s) ds \int_{-s-\sigma}^{-s-\sigma_0} \frac{1}{2}(-s-\sigma_0-y)^2 u_{xx}(y,s) dy,
\end{array}
\]
so that the we have 
\[
|H(t)|\le\frac{1}{2}\left(\int_{0}^{t}(t-s)ds\int_{-s-\sigma}^{-s-\sigma_0}(-s-\sigma_0-y)|u_{xx}(y,s)|^pdy\right)^{1/p} I(t)^{1-1/p},
\]
where
\[
\begin{array}{rl}
I(t):=
&\d\int_{0}^{t} (t-s) ds \int_{-s-\sigma}^{-s-\sigma_0} (-s-\sigma_0-y)^{(2p-1)/(p-1)}dy\\
=&\d\frac{p-1}{2(3p-2)} (\sigma-\sigma_0)^{(3p-2)/(p-1)}t^2\\
=&\d \frac{p-1}{3p-2} 2^{(2p-1)/(p-1)}\sigma_1^{(3p-2)/(p-1)} t^2.
\end{array}
\]
Therefore, combining (\ref{H''bound}) and the estimate of $|H(t)|$ above, we obtain that
\begin{equation}
\label{ineq:H}
H''(t)\ge D t^{1-2p} H(t)^p\quad\mbox{for}\ t\ge\sigma_1,
\end{equation}
where
\[
D:=2^{-p}\left(\frac{p-1}{3p-2}\right)^{1-p}\sigma_1^{-(3p-2)}>0.
\]

\par
From now on, we employ a different method to investigate (\ref{est:H}) and (\ref{ineq:H})
from Rammaha \cite{Rammaha95, Rammaha97},
which is a combination of ordinary differential inequality and an iteration argument
firstly introduced in Lai and Takamura \cite{LT18}.
Namely, assume that
\begin{equation}
\label{j-step}
H(t) \geq C_j(t-\sigma_1)^{a_j}t^{(1-2p)b_j}\quad\mbox{for}\ t\ge\sigma_1,
\end{equation}
where $C_j>0$, $a_j\ge0$, $b_j\ge0$ $(j\in\N)$.
Our first step (\ref{est:H}) should be the case of $j=1$ as $C_1=B,\ a_1=2,\ b_1=0$.
It follows from (\ref{ineq:H}) and $H'(0)=H''(0)=0$ that
\begin{equation}
\label{frame}
H(t)\ge D\int_{\sigma_1}^{t}ds\int_{\sigma_1}^{s} \tau^{1-2p} H(\tau)^p d\tau
\quad\mbox{for}\ t\ge\sigma_1.
\end{equation}
This is the iteration frame.
Substituting $H$ in the right-hand side of (\ref{frame}) by (\ref{j-step}), we have that
\[
\begin{array}{ll}
H(t)
&\d\ge DC_j^p \int_{\sigma_1}^{t} ds \int_{\sigma_1}^{s} \tau^{(1-2p)(pb_j+1)} (\tau-\sigma_1)^{pa_j} d\tau\\
&\d\ge DC_j^p t^{(1-2p)(pb_j+1)} \int_{\sigma_1}^{t} ds \int_{\sigma_1}^{s} (\tau-\sigma_1)^{pa_j} d\tau\\
&\d= \frac{DC_j^p}{(pa_j+1)(pa_j+2)} (t-\sigma_1)^{pa_j+2}t^{(1-2p)(pb_j+1)}\\
&\d\ge\frac{DC_j^p}{(pa_j+2)^2} (t-\sigma_1)^{pa_j+2}t^{(1-2p)(pb_j+1)}.
\end{array}
\]
Hence (\ref{j-step}) is true for all $j\in\N$ if the sequence $\{C_j\}, \{a_j\}, \{b_j\}\ (j\in \N)$ satisfy
\[
\left\{
\begin{array}{ll}
\d C_{j+1}\ge \frac{DC_j^p}{(pa_j+2)^2}, & C_1 = B,\\
a_{j+1}=pa_j + 2, & a_1 = 2,\\
b_{j+1}=pb_j + 1, & b_1 = 0.
\end{array}
\right.
\]
Solving $a_j$ and $b_j$, we have
\[
a_j=\frac{2}{p-1}(p^{j} - 1),\quad
b_j= \frac{1}{p-1}(p^{j-1} - 1)
\]
for $j\in\N$, so that $C_j$ satisfies
\[
C_{j+1}\ge\frac{DC_j^p}{a_{j+1}^2}
\ge\frac{(p-1)^2DC_j^p}{4p^{2(j+1)}}
\]
which implies that
\[
\begin{array}{ll}
\log C_{j+1}
&\d\ge p\log C_j-\log p^{2(j+1)} + \log \frac{D(p-1)^2}{4}\\
&\d\ge p(p\log C_{j-1}-\log p^{2j} + \log\frac{D(p-1)^2}{4}\\
&\d\quad + p\log C_j-\log p^{2(j+1)} + \log \frac{D(p-1)^2}{4}\\
&\d\ge p^j\log C_1-\sum_{i = 1}^{j}  p^{j-i}\left(\log p^{2(i+1)}-\log\frac{D(p-1)^2}{4}\right).
\end{array}
\]
Therefore $C_j$ should be defined by
\[
C_j:=\exp\{p^{j-1}(\log B-E)\},
\]
where
\[
E:=\sum_{i = 1}^{\infty}\frac{1}{p^i}\left(\log p^{2i}+\left|\log\frac{D(p-1)^2}{4}\right|\right)<\infty.
\]
We note that the convergence of the infinite sum can be verified by d'Alembert criterion. 

\par
Turning back to (\ref{j-step}), we have that
\begin{equation}
\label{est:final}
H(t)\ge\frac{t^{(2p-1)/(p-1)}}{(t-\sigma_1)^{2/(p-1)}}\exp\{p^{j-1}J(t)\}
\quad\mbox{for}\ t\ge\sigma_1,\ j\in\N,
\end{equation}
where $J$ is defined by
\[
J(t):= \log(t-\sigma_1)^{(2p)/(p-1)} - \log t^{(2p-1)/(p-1)}+\log B-E.
\]
By virtue of (\ref{est:final}),
we reach to a contradiction $H(t)=\infty$ under the assumption that $u$ is a classical solution
by taking limit as $j\rightarrow\infty$
if $t\in[\sigma_1,T)$ satisfies
\begin{equation}
\label{condition}
\frac{B(t-\sigma_1)^{(2p)/(p-1)}}{t^{(2p-1)/(p-1)}\exp E} > 1
\end{equation}
which yields $J(t)>0$.

\par
Now, assume further that
\[
t\ge2\sigma_1.
\]
Then (\ref{condition}) follows from
\[
\frac{Bt^{1/(p-1)}}{2^{(2p)/(p-1)}\exp E} > 1
\]
which is equivalent to
\[
t>2^{3p-1}\left(\frac{\exp E}{C_f}\right)^{p-1} \e^{-(p-1)}.
\]
Therefore, if $\e_0$ is defined to satisfy
\[
2^{3p-1}\left(\frac{\exp E}{C_f}\right)^{p-1} \e_0^{-(p-1)} = 2\sigma_1,
\]
then the classical solution of $u$ as well as $v$ of (\ref{IVP}) cannot exists
as far as $T$ satisfies
\begin{equation}
\label{conclusion}
T>2^{3p-1}\left(\frac{\exp E}{C_f}\right)^{p-1} \e^{-(p-1)},
\end{equation}
where $\e$ is any positive number of $0<\e\le\e_0$.

\par
Finally, we shall prove the theorem under the assumption (\ref{positive2}) instead of (\ref{positive1}).
In this case, viewing (\ref{H''}), we have that
\[
H''(t)\ge\frac{\e}{2}\int_{-t-\sigma}^{-t-\sigma_0}dx\int_{x-t}^{x+t}dy\int_{-\sigma}^yg(z)dz
\quad\mbox{for}\ t\ge0.
\]
We note that there is a freedom to fix $\sigma_0\in(0,\sigma)$ in this case
as we don't assume (\ref{C_f}).
Here we may assume that
\[
\sigma>3\sigma_0\quad\mbox{and}\quad\int_{-\sigma}^{-(\sigma+\sigma_0)/2}g(z)dz>0.
\]
Then it follows from
\[
x-t\le-2t-\sigma_0
\]
that
\[
H''(t)\ge\frac{\e}{2}\int_{-t-(\sigma+\sigma_0)/2}^{-t-\sigma_0}dx
\int_{-(\sigma+\sigma_0)/2}^{x+t}dy\int_{-\sigma}^yg(z)dz
\quad\mbox{for}\ t\ge\frac{\sigma_1}{2},
\]
where $\sigma_1$ is the one in (\ref{sigma_1}).
Because
\[
-2t-\sigma_0\le-\frac{\sigma+\sigma_0}{2}
\]
is equivalent to $t\ge\sigma_1/2$.
We also note that $\sigma_0<\sigma_1<\sigma$ holds.
Hence we have that
\[
H''(t)\ge\frac{\e}{2}\int_{-t-(\sigma+\sigma_0)/2}^{-t-\sigma_0}dx
\int_{-(\sigma+\sigma_0)/2}^{x+t}dy\int_{-\sigma}^{-(\sigma+\sigma_0)/2}g(z)dz
\quad\mbox{for}\ t\ge\frac{\sigma_1}{2}
\]
which implies
\[
H''(t)\ge C_g\e\quad\mbox{for}\ t\ge\frac{\sigma_1}{2},
\]
where
\[
C_g:=\frac{\sigma_1^2}{4}\int_{-\sigma}^{-(\sigma+\sigma_0)/2}g(z)dz>0.
\]
Since $H''(t)\ge0$ for $t\ge0$ and $H(0)=H'(0)=0$,
we have $H(\sigma_1/2),H'(\sigma_1/2)\ge0$ which yields that
\[
H(t)\ge\frac{C_g}{2}\e\left(t-\frac{\sigma_1}{2}\right)^2\quad\mbox{for}\ t\ge\frac{\sigma_1}{2}.
\]
Therefore we obtain the first step of the iteration,
\[
H(t)\ge\frac{C_g}{8}\e t^2\quad\mbox{for}\ t\ge\sigma_1
\]
instead of (\ref{est:H}).
The conclusion can be established by completely same way as in the previous case
in which $C_f$ is replaced with $C_g/8$.
The proof is now completed.
\hfill$\Box$

%%%%%%%%%%%%%%%%%%%%%%%%%%%%%%%%%%%%%%%%%%%%%%%%%%
%%%%%%%%%%%%%%%%%%%%% SECTION %%%%%%%%%%%%%%%%%%%%%%%
%%%%%%%%%%%%%%%%%%%%%%%%%%%%%%%%%%%%%%%%%%%%%%%%%%

\section{Proof of Theorem \ref{thm:t-derivative}}
The basic strategy in this section is almost the same as John \cite{John81, John_book}.
We employ the contradiction argument again
and, similarly to the previous section, we may assume $A=p$.
The finite propagation speed of the wave (\ref{supp_v}) is still valid,
but in this case we define a new function $u$ by
\begin{equation}
\label{u2}
u(x,t):=\int_{0}^{t}v(x,s)ds.
\end{equation}
Then the equation, $v_{tt}-v_{xx}=p|v_t|^{p-2}v_tv_{tt}$, yields that
\[
\frac{\p}{\p t}\left(u_{tt} - u_{xx} - |u_{tt}|^p\right)=0
\quad\mbox{in}\ \R \times (0,T).
\]
Since it follows from (\ref{supp_v}) in which $t$ is replaced with $s\in[0,t]$ that
\begin{equation}
\label{supp_u2}
\mbox{supp}\ u(x,t)\subset\{(x,t)\in\R\times[0,T):|x|\le t+\sigma\},
\end{equation}
we have, by integration in time, a new problem;
\begin{equation}
\label{eq:u2}
\left\{
\begin{array}{ll}
u_{tt} - u_{xx} = |u_{tt}|^p+\e g-\e^pg^p & \mbox{in}\ \R \times (0,T),\\
u(x,0) = 0,\ u_t(x,0)=\e f(x), & x\in\R
\end{array}
\right.
\end{equation}
because of
\[
u_{tt}(x,0)=v_t(x,0)=\e g(x),\quad u_{xx}(x,0)=\int_0^0v_{xx}(x,s)ds=0.
\]
Hence, similarly to (\ref{eq:integral}), $u$ has to satisfy
\begin{equation}
\label{eq:integral_t}
u=u_0+L(|u_{tt}|^p+\e g-\e^pg^p)\quad\mbox{in}\ \R\times[0,T),
\end{equation}
where
\[
u_0(x,t)=\frac{\e}{2}\int_{x-t}^{x+t}f(y)dy\ge0.
\]
In this section, we assume that
\begin{equation}
\label{smallness}
\e^{p-1}\|g\|_{L^\infty(\R)}^{p-1}\le\frac{1}{2}
\end{equation}
which yields that
\begin{equation}
\label{positive3}
\e g-\e^pg^p=\e g\{1-\e^{p-1}g^{p-1}\}\ge\frac{1}{2}\e g
\quad\mbox{in}\ \R.
\end{equation}
\par
First, we shall prove the theorem under the assumption (\ref{positive1}).
From now on, we set
\begin{equation}
\label{direction}
t=x+\sigma\quad\mbox{for}\ x\ge\sigma.
\end{equation}
Then, we have that
\[
u_0(x,t)=\frac{\e}{2} \int_{-\sigma}^{\sigma}f(y)dy=:\wt{C_f}\e>0
\]
and
\[
\begin{array}{ll}
L(|u_{tt}|^p)(x,t)
&\d= \frac{1}{2} \int_{0}^{t} ds \int_{x-t+s}^{x+t-s} |u_{tt}(y,s)|^p dy\\
&\d\ge\frac{1}{2} \int_{\sigma}^{x} dy \int_{y-\sigma}^{y+\sigma} |u_{tt}(y,s)|^p ds.
\end{array}
\]
We shall estimate the integrand of the right-hand side in the inequality above.
(\ref{supp_u2}) yields that
\[
u(y,y+\sigma)
=u(y,y-\sigma)+\int_{y-\sigma}^{y+\sigma}u_t(y,s)ds
= \int_{y-\sigma}^{y+\sigma} u_t(y,s)ds,
\]
so that we have
\[
u(y,y+\sigma)= \int_{y-\sigma}^{y+\sigma} (y+\sigma-s) u_{tt}(y,s)ds
\quad\mbox{for}\ y\ge\sigma
\]
by integration by parts.
Hence it follows from H\"older's inequality that
\[
\begin{array}{ll}
 |u(y,y+\sigma)|
 &\d\le\int_{y-\sigma}^{y+\sigma} (y+\sigma-s) |u_{tt}(y,s)|ds\\
 &\d\le\left(\int_{y-\sigma}^{y+\sigma} (y+\sigma-s)^{p/(p-1)} ds\right)^{1-1/p}
\left( \int_{y-\sigma}^{y+\sigma}  |u_{tt}(y,s)|^p ds\right)^{1/p}.
\end{array}
\]

\par
Therefore, setting
\[
U(x):=u(x,x+\sigma)\quad\mbox{for}\ x\ge\sigma,
\]
we obtain by (\ref{eq:integral_t}) and (\ref{positive3}) that
\begin{equation}
\label{ineq:U}
U(x) \ge\wt{C_f}\e+ \frac{1}{2F} \int_{\sigma}^{x} U(y)^pdy
\quad\mbox{for}\ x\ge\sigma,
\end{equation}
where a constant $F$ is defined by
\[
F:=\left(\frac{p-1}{2p-1}(2\sigma)^{(2p-1)/(p-1)}\right)^{p-1}>0.
\]
If $W$ is a solution of
\[
W(x)=\wt{C_f}\e+\frac{1}{2F}\int_{\sigma}^{x}W(y)^pdy,
\]
Gronwall's inequality yields that
\begin{equation}
\label{comparison}
U(x)\ge W(x)\quad\mbox{for}\ x\ge\sigma.
\end{equation}
Since $W$ has an expricite form,
\[
W(x)=\left\{\frac{2F(\wt{C_f})^{p-1}\e^{p-1}}{2F-(p-1)(\wt{C_f})^{p-1}\e^{p-1}(x-\sigma)}\right\}^{1/(p-1)},
\]
our setting (\ref{direction}) says that $t$ has to satisfy
\[
(2\sigma\le)\ t<\frac{2F(\wt{C_f})^{-(p-1)}}{p-1}\e^{-(p-1)}+2\sigma
\]
whenever $U$ exists.
It means that, if $\e_1$ is defined to satisfy
\[
\frac{2F(\wt{C_f})^{-(p-1)}}{p-1}\e_1^{-(p-1)}=\sigma
\]
in this case, the classical solution of $u$ as well as $v$ cannot exists
as far as $T$ satisfies
\[
T >\frac{6F(\wt{C_f})^{-(p-1)}}{p-1}\e^{-(p-1)}\quad \mbox{for}\ 0 < \e\le\e_1.
\]
In view of (\ref{smallness}), the theorem in this case is established by taking
\[
\e_0=\min\{\e_1,2^{-1/(p-1)}\|g\|_{L^\infty(\R)}^{-1}\}.
\]
we note that $g\equiv0$ is included in this case.

\par
Finally, we shall prove the theorem under the assumption (\ref{positive2}) instead of (\ref{positive1}).
It follows from (\ref{positive3}) that
\[
L(\e g-\e^pg^p)(x,t)\ge\frac{\e}{4}\int_0^tds\int_{x-t+s}^{x+t-s}g(y)dy.
\]
Then, setting (\ref{direction}) again, we have that
\[
\begin{array}{ll}
L(\e g-\e^pg^p)(x,t)
&\d\ge\frac{\e}{4}\int_{-\sigma}^xg(y)dy\int_0^{y+\sigma}ds\\
&\d\quad+\frac{\e}{4}\int_x^{x+2\sigma}g(y)dy\int_0^{-y+x+2\sigma}ds\\
&\d\ge\frac{\e}{4}\int_{-\sigma}^x(y+\sigma)g(y)dy.
\end{array}
\]
Hence, neglecting $u^0$ in (\ref{eq:integral_t}), we obtain that
\[
U(x)\ge\wt{C_g}\e+\frac{1}{2F}\int_\sigma^xU(y)^pdy
\quad\mbox{for}\ x\ge\sigma
\]
instead of (\ref{ineq:U}),
where
\[
\wt{C_g}:=\frac{1}{4}\int_{-\sigma}^\sigma(y+\sigma)g(y)dy>0.
\]
Therefore the same conclusion holds,
in which $\wt{C_f}$ is replaced with $\wt{C_g}$. 
The proof is now completed.
\hfill$\Box$

%%%%%%%%%%%%%i%%%%%%%%%%%%%%%%%%%%%%%
%%%%%%%%%%%% Concluding remark %%%%%%%%%%%%%%%
%%%%%%%%%%%%%%%%%%%%%%%%%%%%%%%%%%%%%
\section{Concluding remarks and future problems}
\par
It is clear that the optimality of Theorem \ref{thm:main} is open for general $p$
except for even integers.
More than this, we are interested in essentially different type of quasilinear wave equations
from (\ref{IVP}), for example, $u_{tt}-u_{xx}=u_t^{p-1}u_{xx}$ with integer $p(\ge2)$,
to which the general theory can be applied.
Because, as we see, the method of the proof of the blow-up theorem
for the nonlinear term of  $|u_x|^p$ is different from the one of $|u_t|^p$.
For this kind of the mixed case, we don't know what will happen on the blow-up part till now.

\par
Another interest may go to higher dimensional case.
For $n$ dimensional version of semilinear equation (\ref{eq:STT}), say
\[
v_{tt}-\Delta_x v=|\nabla_xv|^p\quad\mbox{in}\ \R^n\times(0,T),\ n\ge2,
\]
Shao, Takamura and Wang \cite{STW} have recently studied the lifespan estimates
and obtained that
\[
T(\e)\le
\left\{
\begin{array}{ll}
C\e^{2(p-1)/\{2-(n-1)(p-1)\}} & \mbox{for}\ 1<p<p_G(n),\\
\exp\left(C\e^{-(p-1)}\right) & \mbox{for}\ p=p_G(n),
\end{array}
\right.
\]
where $C$ is a positive constant independent of $\e$,
and $p_G(n):=(n+1)/(n-1)$ is the so-called Glassey exponent.
Note that the upper bound in (\ref{est:STT}) coincides with the first line of the estimate above,
the sub-critical case, when one formally set $n=1$.
The original Glassey conjecture is organized for the equation (\ref{eq:Zhou}) in $n$ space dimensions.
See the introduction of \cite{STW} for its whole histories and references therein.
One may expect that the same lifespan estimates can be established
even for some kind of high dimensional version of (\ref{IVP}), for example,
\[
v_{tt}-\frac{n-1}{r}v_r-v_{rr}=|v_r|^{p-2}v_rv_{rr} \quad\mbox{in}\ (0,\infty)\times(0,T),\ n\ge2,
\]
for the radial unknowns $v(x,t)=v(r,t),\ r:=|x|$.

%%%%%%%%%%%%%i%%%%%%%%%%%%%%%%%%%%%%%
%%%%%%%%%%%% Acknowledgement %%%%%%%%%%%%%%%
%%%%%%%%%%%%%%%%%%%%%%%%%%%%%%%%%%%%%
\section*{Acknowledgement}
\par
The second author is partially supported
by the Grant-in-Aid for Scientific Research (A) (No. 22H00097), 
Japan Society for the Promotion of Science.
The authors appreciate the reviewers for their many helpful comments
which make the manuscript completed.

%%%%%%%%%%%%%%%%%%%%%%%%%%%%%%%%%%%%%%
%%%%%%%%%%%% References %%%%%%%%%%%%%%%%%%%%
%%%%%%%%%%%%%%%%%%%%%%%%%%%%%%%%%%%%%%

\bibliographystyle{plain}

\begin{thebibliography}{20}

\bibitem{John81}{F. John},
{\it Blow-up for quasilinear wave equations in three space dimensions},
Comm. Pure Appl. Math. {\bf 34} (1981), no. 1, 29-51.

\bibitem{John_book}{F. John},
\lq\lq Nonlinear Wave Equations, Formation of Singularities",
ULS Pitcher Lectures in Mathematical Science, Lehigh University,
American Mathematical Society, Providence, RI, 1990.

\bibitem{KSTT}{R. Kido, T. Sasaki, S. Takamatsu and H. Takamura},
{\it The generalized combined effect for one dimensional wave equations
with semilinear terms including product type},
J. Differential Equations, {\bf 403} (2024), 576-618.

\bibitem{STW}{K. Shao, H. Takamura and C. Wang},
{\it Blow-up of solutions to semilinear wave equations with spatial derivatives},
arXiv: 2406.02098, to appear in Discrete and Continuous Dynamical Systems.

\bibitem{KMT23}{S. Kitamura, K. Morisawa and H. Takamura},
{\it Semilinear wave equations of derivative type
with spatial weights in one space dimension},
Nonlinear Analysis, RWA. {\bf 72} (2023), Paper No. 103764.

\bibitem{LT18}{N.-A. Lai and H. Takamura},
{\it Blow-up for semilinear damped wave equations with subcritical exponent in the scattering case},
Nonlinear Anal. 168 (2018), 222-237. 

\bibitem{LYZ91}{T.-T. Li (D.-Q. Li), X. Yu and Y. Zhou},
{\it Dur\'ee de vie des solutions r\'eguli\`eres pour les \'equations des ondes non lin\'eaires unidimensionnelles} (French),
C. R. Acad. Sci. Paris S\'er. I Math., {\bf 312} (1991), no. 1, 103-105.

\bibitem{LYZ92}{T.-T. Li, X. Yu and Y. Zhou},
{\it Life-span of classical solutions to one-dimensional nonlinear wave equations},
Chinese Ann. Math., Ser. B, {\bf13} (1992), no. 3, 266-279. 

\bibitem{MST23}{K. Morisawa, T. Sasaki and H. Takamura},
{\it The combined effect in one space dimension beyond the general theory
for nonlinear wave equations},
Commun. Pure Appl. Anal. {\bf 22} (2023), no. 5, 1629-1658.

\bibitem{MST23e}{K. Morisawa, T. Sasaki and H. Takamura},
{\it Erratum to \lq\lq The combined effect in one space dimension
beyond the general theory for nonlinear wave equations"},
Commun. Pure Appl. Anal. {\bf 22} (2023), no. 10, 3200-3202. 

\bibitem{Rammaha95}{M. A. Rammaha},
{\it Upper bounds for the life span of solutions to systems of nonlinear wave equations
in two and three space dimensions},
Nonlinear Anal. {\bf 25} (1995), no. 6, 639-654.

\bibitem{Rammaha97}{M. A. Rammaha},
{\it A note on a nonlinear wave equation in two and three space dimensions},
Comm. Partial Differential Equations {\bf 22} (1997), no. 5-6, 799-810. 

\bibitem{Sasaki18}{T. Sasaki},
{\it Regularity and singularity of the blow-up curve
for a wave equation with a derivative nonlinearity},
Adv. Differential Equations {\bf 23} (2018), no. 5-6, 373408.

\bibitem{STT23}{T. Sasaki, S. Takamatsu, H. Takamura},
{\it The lifespan of classical solutions of one dimensional wave equations
with semilinear terms of the spatial derivative},
AIMS Math. {\bf 8} (2023), no. 11, 25477-25486.

\bibitem{Sugiyama2022}{Y. Sugiyama},
{\it Formation of singularities for a family of 1D quasilinear wave equations},
Indiana Univ. Math. J. {\bf 71} (2022), no. 6, 2529-2549.

\bibitem{Takamatsu}{S. Takamatsu},
{\it Improvement of the general theory for one dimensional nonlinear wave equations
related to the combined effect},
arXiv:2308.02174.

\bibitem{Takamura}{H. Takamura},
{\it Recent developments on the lifespan estimate for classical solutions of nonlinear wave equations
in one space dimension}, arXiv:2309.08843, to appear in Advanced Studies in Pure Mathematics, MSJ.

\bibitem{Zhou01}{Y. Zhou},
{\it Blow up of solutions to the Cauchy problem for nonlinear wave equations},
Chinese Ann. Math. Ser. B, {\bf22} (2001), no. 3, 275-280.


\end{thebibliography}

\end{document}